# The Two Ignored Components of Random Variation

Haim Shore

Department of Industrial Engineering and Management, Faculty of Engineering Sciences, Ben-Gurion University of the Negev, POB 653, Beer-Sheva 84105, Israel. shor@bgu.ac.il.

**July, 2015**

## ABSTRACT

**Short Abstract:** A random phenomenon may have two sources of random variation: an unstable "identity" and a set of external variation-generating factors. When only a single source is active, two mutually exclusive extreme scenarios may ensue that result in the exponential or the normal, the only truly univariate distributions. All other supposedly univariate random variation observed in nature is truly bivariate. In this article, we elaborate on this new paradigm for random variation and develop a general bivariate distribution to reflect it. It is shown that numerous current univariate distributions are special cases of an approximation to the new bivariate distribution.

**Extended Abstract:** We first show that the exponential and the normal are special cases of a single distribution represented by a Response Modeling Methodology model. We then develop a general bivariate distribution commensurate with the new paradigm, its properties are discussed and its moments developed. An approximating assumption results in a *univariate* general distribution that is shown to include as exact special cases widely used distributions like generalized gamma, log-normal, F, t, and Cauchy. Compound distributions and their relationship to the new paradigm are addressed. Empirical observations that comply with predictions derived from the new paradigm corroborate its scientific validity.

## 1. INTRODUCTION

In a recent article (Shore, 2015), the need for unification of models of random variation within the science of statistics had been addressed. It was emphasized that this need and unification efforts that may ensue are not unlike similar efforts exercised throughout all branches of modern science in the last century or so. An example is modern-day physics, where attempts at unifying the fundamental forces of nature have been ongoing for over a century, resulting most recently in the development of super-string theory.

At the afore-cited paper, a certain paradigm had been introduced for developing a general model of random variation and an example for such a model developed and extensively explored. In this article, we introduce a new paradigm for general modeling of random variation, based on a distinction between two independent sources of random variation, heretofore ignored in the statistical literature. In Section 2 we explain the basic paradigm that forms the theoretical basis for the new general model of random variation and develop a new bivariate distribution based on the new paradigm. A general expression for the moments of this distribution is developed. A univariate distribution that approximates the bivariate distribution is derived in Section 3, assuming that the two r.v.s are perfectly linearly correlated. Some current univariate distributions are shown to be special cases of the new univariate approximation. Sections 4 and 5 address compound distributions and Section 6 the gamma distributions, all in relationship to the two-ignored-components principle. A concluding Section 7 summarizes this article.

## 2. THE NEW PARADIGM AND THE ALLIED BIVARIATE DISTRIBUTION

### 2.1 An initiating example

As a departure point for understanding the distinction between the two prime and indistinguishable sources of random variation, consider duration times of jobs carried out in a typical car garage. Broadly speaking, these may be divided into two categories of jobs:

* Routine repetitive tasks, like regular maintenance carried out at fixed time intervals as specified by the car producer;

* Repair jobs carried out on damaged cars (whatever the source of the damage may be like cars collision, careless driving or a sudden part failure).

Constructing histograms of duration times of jobs, sampled from the two categories, would reveal a surprising distinction: histograms of repetitive maintenance jobs tend to look symmetric; jobs performed on damaged cars tend to look as though sampled from the exponential distribution.

Why is that?

The explanation is simple: repetitive jobs have typical work contents, emanating from highly-detailed specification of a series of maintenance procedures, required and mandated by the producer. In other words, the random variable (henceforth r.v.), denoted "Job Duration", is random not due to "unstable" work contents (absent "identity") but rather due to numerous small "outside" factors, none of which exercise decisive effect on the final value of a realization of the random phenomenon, namely, on job duration. Examples for these factors are varying work rates by different garage employees, unexpected interruptions or delays while performing the job (for example, unavailability of a required spare part) and so on. All these small effects combined will expectedly produce a normal distribution for the job duration. However, all jobs in this category share common work contents, in other words, they share a common "identity". No memory is lost in transition from one job to another.

Conversely, observe a typical sample of jobs carried out on damaged cars (repair jobs). A different picture is revealed: Each job has its own unique series of required work elements, rendering jobs belonging to this category devoid of any typical "identity". There is no "memory" carried over from one job to the next. In other words, job duration is memoryless. This would result in sample-based histograms resembling the pdf of the memoryless exponential distribution.

In real life the two "pure" extreme scenarios, just delineated, rarely exist. Most random phenomena observed in nature are outcomes of varying degrees of mixture of the two "pure" scenarios of a "memory-full" identity (normal) and a memoryless identity (exponential) or, more accurately, an identity with partial memory. Therefore, it stands to reason that a general bivariate distribution, integrating the two independent sources of random variation, would probably provide an adequate general model of random variation from which current widely-used non-mixture univariate statistical distributions may be derived as approximations. In the next subsections 2.2 and 2.3, the new bivariate distribution is constructed and its properties explored and demonstrated. Univariate approximations to this distribution will be developed and demonstrated in Section 3.

**2.2 The normal and the exponential are a single distribution represented by an RMM model**

Let U and V be an exponential and a normal r.v.s, with respective parameters {b} and {µ,σ}. One realizes that the probability density function (pdf) of both distributions derive from a single RMM model (Shore, 2005, 2011, 2012):

$$f(z;\alpha,\kappa,\lambda) = (\kappa)\exp\{-(\alpha/\lambda)(z^\lambda - 1)\}, \{-L < z < \infty, L > 0; 1 \leq \lambda \leq 2\}, \quad (1)$$

where $\kappa$ is a normalizing parameter (ensuring that f is a proper pdf), {$\alpha, \lambda, L$} are parameters,

$$Z = \frac{X - c}{\sigma} \quad (2)$$

is a linear transformation of the original-scale r.v., X, c is some central-tendency parameter and $\sigma$ is the standard deviation. It is easy to show that the mode of (1) is obtained at z=0. Provided $\lambda \geq 1$, the pdf at the mode, denote the latter by M, is

$$f(0;\alpha,\kappa,\lambda) = (\kappa)\exp(\alpha/\lambda), \{1 \leq \lambda \leq 2\}. \quad (3)$$

Eq. 3 implies that z in (1) should represent the original r.v., X, standardized around the mode (rather than around the mean). This would ensure preservation of the mode both on the original scale (X) and on the standardized scale (Z). Therefore, z in (1), in terms of the original scale, x, is:

$$z = \frac{x - M}{\sigma}, \qquad (4)$$

with M being the mode of the distribution of X. It is easy to realize that (1) represents both the exponential and the normal distributions:

* For {$\alpha=1$, $\lambda=1$, L= (M/$\sigma$)= 0}, (1) is the exponential pdf;

* For {$\alpha=1$, $\lambda= 2$, L= $\infty$}, (1) is the normal pdf.

With {M=0, $\sigma=1$} in (4), eq. 1 represents the pdf of X (the r.v. expressed in the original scale). In that case, (1) represents additional distributions (with $\lambda$ allowed to approach zero asymptotically):

* For {$\lambda \to 0$, $\alpha>1$}, (1) converges to the Pareto pdf;

* For {$\lambda \to 0$, $\alpha<1$}, (1) converges to the power function pdf;

* For $\alpha=0$, (1) becomes a uniform distribution on the interval {0, 1/$\kappa$}.

The k-th non-central moment (moment around zero) of (1) is:

$$E(Z^k) = (\frac{1}{\lambda}) \kappa e^{\alpha/\lambda} \left( \frac{\alpha(-1)^\lambda}{\lambda} \right)^{-\frac{k+1}{\lambda}} \left( (-1)^{k+1} \Gamma\left( \frac{k+1}{\lambda}, \frac{\alpha(-L)^\lambda}{\lambda} \right) + \left( \left((-1)^\lambda\right)^{\frac{k+1}{\lambda}} + (-1)^k \right) \Gamma\left( \frac{k+1}{\lambda} \right) \right)$$

(5)

where $\Gamma(\xi)$ is the Euler gamma function and $\Gamma(\gamma,\xi)$ is the incomplete gamma function. Note that (5) is function of four parameters, {$\alpha,\kappa,\lambda,L$}, with k being the order of the moment calculated. Putting k=0 and equating: $E(Z^0)=1$, parameter $\kappa$ can be expressed as function of the other parameters:

$$\kappa = \frac{\lambda e^{-\frac{\alpha}{\lambda}} \left( \frac{\alpha(-1)^\lambda}{\lambda} \right)^{1/\lambda}}{\left( \left((-1)^\lambda\right)^{1/\lambda} + 1 \right) \Gamma\left( \frac{1}{\lambda} \right) - \Gamma\left( \frac{1}{\lambda}, \frac{\alpha(-L)^\lambda}{\lambda} \right)} \qquad (6)$$

Inserting (6) into (5), the latter becomes, for a given k, function of three parameters only, namely, $\{\alpha, \lambda, L\}$.

As earlier expounded, the two ignored components of random variation, at play in any observed steady-state random phenomenon, are tightly linked to the exponential and the normal distributions. In fact, the two are end points of a continuous spectrum of distributions that span different mixtures of the basic "two ignored components", inherent in all observed realizations of random variation. In other words:

* One cannot conceive of an r.v. that does not, at least partially, reflect variation due to the sum total of external numerous small effects, namely, variation represented by a normal component; An exception to this rule is an r.v. pursuing the exponential distribution;

* Conversely, it is difficult to conceive of an r.v. that has no "memory loss" (unstable identity), observable on transition from one realization to another.

There are only two "pure cases", where the two *extreme and rare* scenarios of either complete memory loss, or no memory loss at all, do exist. These are the exponential and the normal cases, respectively. The new paradigm, suggesting a general distribution based on the underlying principle of the "two ignored components", should adequately mirror this reality.

## 2.3 A general model integrating the two-ignored-components of random variation

The normal and the exponential belong to the exponential family of distributions, the pdf of which may be found in any source dealing with generalized linear models (for example, Meyers *et al.*, 2002). Also, mixtures of normal and exponential distributions are well known and have been extensively investigated. Therefore, a natural tendency would deliver representation to the distribution of a r.v. comprising the two-ignored-components, in the form a simple mixture, for example:

$$f(x; w) = \omega * f_U + (1-\omega) * f_V, 0 \leq w \leq 1 ,  \qquad (7)$$

where {U,V} are independent exponential and normal variables, respectively, $f_U$ and $f_V$ are the respective pdfs, and $\omega$ is a weighting coefficient. However, this mixture violates the very nature of the-two-ignored components since it assumes that different mixtures contain weighted sums of independent exponential and normal variables (with varying contributions to the overall weighted average). In reality, once "memory-loss" is not comprehensive but only partial, U ceases to be exponential and concurrently a normal component pops up that starts contributing to the overall random variation. This is due to the fact that the exponential represents total loss of memory, in which case no normal component can be present. Conversely, once "memoryless-ness" is not comprehensive, a normal component is introduced, on the one hand, and on the other hand the component of variation reflecting the degree of stability (or lack thereof) of the "inner identity" ceases to be exponential. In other words, no symmetry exists in the roles played by the two distributions. Therefore a different approach needs to be pursued.

In conformance with the above description of the nature of the two-ignored-components, consider the following model with external effects expressed by a normal multiplicative error, $\varepsilon$:

$$W = U(1+\varepsilon) = (U)(V), U \geq 0, W \geq 0, 1 \leq \lambda \leq 2, \qquad (8)$$

where {U,V} are two independent r.v.s, $\varepsilon$ is normally distributed with zero mean and standard deviation $\sigma_\varepsilon$, and U ($\geq 0$) has pdf given by (1). We also assume W≥0 so that: $\varepsilon \geq -1$, or, equivalently: V≥0.

Consider the joint pdf of {U,V}, based on (1) and taking account of the statistical independence of {U,V} (for simplicity the standardized U and V are denoted $Z_1$ and $Z_2$, respectively):

$$f(z_1, z_2; \lambda_1, \lambda_2, \kappa) = (\kappa)\exp\{-(\alpha_1)(\frac{1}{\lambda_1})[(z_1)^{\lambda_1} - 1] - (\alpha_2)(\frac{1}{\lambda_2})[(z_2)^{\lambda_2} - 1]\}; \alpha_1 + \alpha_2 = 1,$$

(9)

where $\alpha_1$ and $\alpha_2$ are functions of $\lambda_1$ and/or $\lambda_2$, $\kappa$ is a normalizing coefficient and:

$$z_1 = \frac{u - M_1}{\sigma_1} \; ; \; z_2 = \frac{v - M_2}{\sigma_2}, \tag{10}$$

with $\{M_1, \sigma_1\}$ and $\{M_2, \sigma_2\}$ being parameters associated with $Z_1$ and $Z_2$, respectively. A special case of (9), representing the joint distribution of U and V in (8), is:

$$f(z_1, z_2; \kappa, \lambda) = (\kappa) \exp\{-(2-\lambda)(\frac{1}{\lambda})[(z_1)^\lambda - 1] - (\lambda - 1)(\frac{1}{2})[z_2^2 - 1]\},$$
$$\{ For \lambda \neq 2 : -\frac{M_1}{\sigma_1} \leq z_1 < \infty \; ; \; -\frac{1}{\sigma_2} \leq z_2 \leq \infty \; ; \; 1 \leq \lambda \leq 2\} \tag{11}$$

where $\{z_1, z_2\}$ are as defined in (10), namely, $\{M_1, \sigma_1\}$ are parameters associated with $Z_1$, the normalized U, $\{M_2, \sigma_2\}$ are parameters associated with the normal $Z_2$ (with $M_2 = E(V) = 1$), the weighting coefficients of (9) in (11) are $\{\alpha_1 = 2-\lambda \; ; \; \alpha_2 = \lambda - 1\}$, $\lambda_2$ of (9) is equal 2 in (11) and $\kappa$ is a normalizing coefficient. As both W and U in (8) are assumed to be non-negative, we expect $\sigma_2$ to be smaller than, say, about 1/6.

Since $\alpha_1 + \alpha_2 = 1$ (as required in (9)), (11) implies that the two coefficients concurrently fulfill two functions:

- Invoking the exponential and the normal distributions as endpoints of a continuous spectrum extending from total lack of memory (exponential: $\lambda = 1$) to total preservation of memory (normal: $\lambda = 2$);
- Representing the combined effect of the two-ignored-components as components in a single *weighted* average (observe (11)).

The mean of W, from (8), is:

$E(W) = \mu_W = E(U)E(V) = E(U) = \mu_1$ (since $E(V) = 1$) , (12)

and the variance of W is:

$$Var(W) = E(W^2) - [E(W)]^2 = E(U^2)E(V^2) - [E(U)]^2 . \tag{13}$$

From (4) we have for Z:

$$E(Z) = E[\frac{(X-\mu)+(\mu-M)}{\sigma}] = \frac{(\mu-M)}{\sigma}$$

$$E(Z^2) = E\{[\frac{(X-\mu)+(\mu-M)}{\sigma}]^2\} = 1 + (\frac{\mu-M}{\sigma})^2 \qquad (14)$$

The standardized U and V, namely, $Z_1$ and $Z_2$, have pdf with parameters (refer to (9)):

For $Z_1$: $\{\alpha_1, \lambda_1, L_1\} = \{2-\lambda, \lambda, M_1/\sigma_1\}$;

For $Z_2$: $\{\alpha_2, \lambda_2, L_2\} = \{\lambda-1, 2, M_2/\sigma_2\} = \{\lambda-1, 2, 1/\sigma_2\}$. $\qquad (15)$

Therefore, from (11)-(15):

$$E(Z_1) = \frac{\mu_1 - M_1}{\sigma_1} =$$

$$(\frac{1}{\lambda})\kappa_1 e^{\alpha_1/\lambda} \left(\frac{\alpha_1(-1)^\lambda}{\lambda}\right)^{-\frac{2}{\lambda}} \left((-1)^2 \Gamma\left(\frac{2}{\lambda}, \frac{\alpha_1(-L_1)^\lambda}{\lambda}\right) + \left(((-1)^\lambda)^{\frac{2}{\lambda}} + (-1)^1\right)\Gamma\left(\frac{2}{\lambda}\right)\right) =$$

$$(\frac{1}{\lambda})\kappa_1 e^{\alpha_1/\lambda} \left(\frac{\alpha_1(-1)^\lambda}{\lambda}\right)^{-\frac{2}{\lambda}} \left(\Gamma\left(\frac{2}{\lambda}, \frac{\alpha_1(-L_1)^\lambda}{\lambda}\right) + \left(-1 + ((-1)^\lambda)^{\frac{2}{\lambda}}\right)\Gamma\left(\frac{2}{\lambda}\right)\right)$$

$$E(Z_1^2) = 1 + (\frac{\mu_1 - M_1}{\sigma_1})^2 =$$

$$(\frac{1}{\lambda})\kappa_1 e^{\alpha_1/\lambda} \left(\frac{\alpha_1(-1)^\lambda}{\lambda}\right)^{-\frac{3}{\lambda}} \left((-1)^3 \Gamma\left(\frac{3}{\lambda}, \frac{\alpha_1(-L_1)^\lambda}{\lambda}\right) + \left(((-1)^\lambda)^{\frac{3}{\lambda}} + (-1)^2\right)\Gamma\left(\frac{3}{\lambda}\right)\right) =$$

$$(\frac{1}{\lambda})\kappa_1 e^{\alpha_1/\lambda} \left(\frac{\alpha_1(-1)^\lambda}{\lambda}\right)^{-\frac{3}{\lambda}} \left(-\Gamma\left(\frac{3}{\lambda}, \frac{\alpha_1(-L_1)^\lambda}{\lambda}\right) + \left(1 + ((-1)^\lambda)^{\frac{3}{\lambda}}\right)\Gamma\left(\frac{3}{\lambda}\right)\right)$$

$$E(Z_2) = \frac{\mu_2 - M_2}{\sigma_2} = 0 =$$

$$(\frac{1}{2})\kappa_2 e^{\alpha_2/2} \left(\frac{\alpha_2(-1)^2}{2}\right)^{-\frac{2}{2}} \left((-1)^2 \Gamma\left(\frac{2}{2}, \frac{\alpha_2(-L_2)^2}{2}\right) + \left(((-1)^2)^{\frac{2}{2}} + (-1)^1\right)\Gamma\left(\frac{2}{2}\right)\right) =$$

$$(\frac{\kappa_2}{\alpha_2})e^{(\alpha_2/2)(1-L_2^2)}$$

$$(16)$$

$$E(Z_2^2) = 1 + (\frac{\mu_2 - M_2}{\sigma_2})^2 = 1 =$$

$$(\frac{1}{2})\kappa_2 e^{\alpha_2/2} \left(\frac{\alpha_2(-1)^2}{2}\right)^{-\frac{3}{2}} \left((-1)^3 \Gamma\left(\frac{3}{2}, \frac{\alpha_2(-L_2)^2}{2}\right) + \left(((-1)^2)^{\frac{3}{2}} + (-1)^2\right) \Gamma\left(\frac{3}{2}\right)\right) =$$

$$\kappa_2 e^{\alpha_2/2} \left(\frac{2}{\alpha_2^3}\right)^{\frac{1}{2}} \left(\sqrt{\pi} - \Gamma\left(\frac{3}{2}, \frac{\alpha_2(L_2)^2}{2}\right)\right)$$

with $\{\kappa_1, \kappa_2\}$ given by (6) ($\kappa$ in (11) is equal $(\kappa_1\kappa_2)$). Note that all equations in (16) are taken from (5) and they are expressed in terms of the parameters given in (15), either implicitly or explicitly. Also note that Var($Z_1$)=Var($Z_2$)=1, an expected result since a change in location of an r.v. does not change its variance.

Parameter $\mu_1$ in (16) is estimated by the *observed* E(W) (relate to (12)). Also from (12) and (13) we may obtain for $\sigma_1^2$:

$$\sigma_1^2 = \frac{Var(W) - [E(W)]^2 \sigma_2^2}{(1+\sigma_2^2)} \qquad (17)$$

(note that $\sigma_2$ is dimensionless by definition; relate to (8))

Eqs. 16 are expressed in terms of the observable mean of W, E(W), the observable variance of W, Var(W), and three unknowns: $\{\lambda, M_1, \sigma_2\}$. The latter may be found by minimizing the sum of squared deviations between the right-hand side and the left-hand side of the four equations of (16).

For any given distribution with known mean, E(W), and known variance, V(W), the resulting fitted distribution would embody the two-ignored-component principle, with estimates obtained by a two-moment-matching procedure.

### 2.4 The distribution of W=U*V

In (3.3) we have derived the joint distribution of $\{Z_1, Z_2\}$, where the latter represent the standardized U and V, respectively (relate to eq. 6). To derive the pdf of W, $f_W(w)$, note that:

$$f_W(w) = \int_{-\infty}^{\infty} f_{u,v}(\frac{w}{v}, v) \frac{1}{|v|} dv, \qquad (18)$$

where $f_{u,v}$ is given by the joint distribution (11). We were unable to derive the pdf of W in (18) in explicit form.

## 3. SIMPLIFIED UNIVARIATE APPROXIMATIONS

### 3.1 A simplified approximate model – $Z_1$ and $Z_2$ are perfectly linearly correlated

Model (8), and the associated joint pdf (11), assume that the *observed* W reflects the concurrently-active and independent $Z_1$ and $Z_2$. Conversely, assume that $Z_1$ and $Z_2$ are perfectly linearly correlated, namely:

$$Z_1 = Z;$$
$$Z_2 = \beta_0 + \beta_1 Z_1 = \beta_0 + \beta_1 Z. \tag{19}$$

This assumption should be considered a simplifying one, leading to a univariate approximation of the true (actual) bivariate scenario, as represented by (8). The approximation allows us to express the pdf of *W* as function of a single variable, z:

$$f(z;\boldsymbol{\theta}) = (\kappa)\exp\{-(\alpha_1)(\frac{1}{\lambda_1})[(z)^{\lambda_1} - 1] - (\alpha_2)(\frac{1}{\lambda_2})[(\beta_0 + \beta_1 z)^{\lambda_2} - 1]\}, \{1 \leq \{\lambda_1, \lambda_2\} \leq 2\}, \tag{20}$$

where $\boldsymbol{\theta}$ is a vector of parameters. It can easily be shown that for (20) to have a single mode we should have: $\beta_0 = 0$. Therefore the two perfectly correlated $\{Z_1, Z_2\}$, as represented in (20), only differ in scale, namely, have different standard deviation. Eq. 20 then becomes:

$$f(z;\boldsymbol{\theta}) = (\kappa)\exp\{-(\alpha_1)(\frac{1}{\lambda_1})[(z)^{\lambda_1} - 1] - (\alpha_2)(\frac{1}{\lambda_2})[(\beta z)^{\lambda_2} - 1]\}, \{1 \leq \{\lambda_1, \lambda_2\} \leq 2\}, \tag{21}$$

with $\beta_1$ replaced by $\beta$. Note, however, that (20) requires $\beta_0 = 0$ to be unimodal only if the constraint on the values of $\{\lambda_1, \lambda_2\}$ is fulfilled. Since (21) is only an approximation to the joint pdf (9), and the latter is based on the exact model (8), we do not expect the two basic constraints of (9), namely,

$$1 \leq \{\lambda_1, \lambda_2\} \leq 2, \quad \alpha_1 + \alpha_2 = 1, \tag{22}$$

to always be preserved in current *univariate* distributions, when the latter are shown to be special cases of (20). However, we expect the basic algebraic structure of (20) to be preserved. Therefore we generalize (20) to obtain:

$$f(z; \boldsymbol{\theta}) = (\kappa) \exp\{-(\alpha_1)(\frac{1}{\lambda_1})[(z)^{\lambda_1} - 1] - (\alpha_2)(\frac{1}{\lambda_2})[(\beta_0 + \beta_1 z)^{\lambda_2} - 1]\}, \{0 \leq \{\lambda_1, \lambda_2\} \leq UL\}, \tag{23}$$

with UL being the upper limit of the parameters $\{\lambda_1, \lambda_2\}$. Again, note that now the constraint on the values of $\{\lambda_1, \lambda_2\}$, as appearing in (21), has been modified.

Here are some examples for the generality of (23) as an approximation to the exact model (9), leading to regarding numerous current univariate distributions as special cases of (23).

Consider the Weibull distribution with pdf:

$$f(w; b, c) = (\kappa)(x/b)^{c-1} \exp[-(x/b)^c], x \geq 0, \tag{24}$$

with $\kappa$ as a normalizing parameter. It is easy to show that (24) is a special case of (23) with $z = (x/b)$ and parameters: $\lambda_1 = 0$; $\lambda_2 = c$; $\alpha_1 = 1-c$; $\alpha_2 = c = 1-\alpha_1$; $\beta_0 = 0$; $\beta_1 = 1$.

Next, consider the generalized gamma distribution with pdf:

$$f(w; \boldsymbol{\theta}) = (\kappa)[(x-a)/(bc^{1/k})]^{kc-1} \exp\{-(kc)[(x-a)/(bc^{1/k})]^k / k\}, x \geq 0, \tag{25}$$

with $\boldsymbol{\theta} = \{a, b, c, k\}$ being a vector of parameters and $\kappa$ is a normalizing parameter. It is easy to show that (25) is a special case of (23) with $z = (x-a)/(bc^{1/k})$ and parameters: $\lambda_1 = 0$; $\lambda_2 = k$; $\alpha_1 = 1-kc$; $\alpha_2 = kc = 1-\alpha_1$; $\beta_0 = 0$; $\beta_1 = 1$. Also, gamma (with k=1, a=0), the exponential (with c=k=1, a=0), Weibull (with c=1, a=0) and the chi-squared variate (a=0, b=2, c=n/2, k=1) are all special cases of this distribution. It is interesting to note that (25), and all its special cases as just enumerated, fulfill the fundamental condition that the two weights, $\{\alpha_1, \alpha_2\}$, in (23) sum up to one.

Next, consider the F distribution:

$$f(w;m,n) = \kappa \frac{w^{m/2-1}}{[1+(\frac{m}{n})w]^{(1/2)(m+n)}}, \tag{26}$$

with κ as a normalizing parameter. It is easy to show that (26) is a special case of (23) with z=w and parameters: $\lambda_1$=0; $\lambda_2$=0; $\alpha_1$=1-(m/2); $\alpha_2$=(m/2)+(n/2); $\beta_0$=1; $\beta_1$=m/n. We realize that this time the second condition in (22) is not fulfilled.

Next consider the log-normal distribution:

$$f(w;\mu,\sigma) = \kappa \exp\{-[\log(w)-\mu]-(1/2)(\frac{\log(w)-\mu}{\sigma})^2], \tag{27}$$

and κ is a normalizing parameter. It is easy to show that (27) is a special case of (23) with z=[log(w)-μ] and parameters: $\lambda_1$=1; $\lambda_2$=2; $\alpha_1$=1; $\alpha_2$=1; $\beta_0$=0; $\beta_1$=1/σ.

Next consider Student's t distribution:

$$f(w;m) = \kappa(1+\frac{w^2}{m})^{-(m+1)/2}, \tag{28}$$

and κ as a normalizing parameter. It is easy to show that (28) is a special case of (23) with z = $w^2$ and parameters: $\alpha_1$=0; $\lambda_2$=0; $\alpha_2$=(m+1)/2; $\beta_0$=1; $\beta_1$=1/m.

Next consider Cauchy distribution:

$$f(w;a,b) = \kappa[1+(\frac{w-a}{b})^2]^{-1}, \tag{29}$$

with κ as a normalizing parameter. It is easy to show that (30) is a special case of (23) with z = $[(w-a)/b]^2$ and parameters: $\alpha_1$=0; $\lambda_2$=0; $\alpha_2$=1; $\beta_0$=1; $\beta_1$=1. We realize that: $\alpha_1+\alpha_2$=1, corroborating once again the validity of (11) as a general model for random variation.

Numerous other distributions also comply with (23). However, we reiterate that all these univariate distributions are considered by us **approximations** to the true bivariate distribution, delivered by pdf (11), since, being univariate, they all fail to comply with the universal two-ignored-components principle of random variation, as expounded in this article. Furthermore, some of the above distributions had been derived as sampling distributions of statistics, which are functions of other statistics assumed to follow the normal distribution. The latter assumption, as claimed here, is a dubious one in light of the two-ignored-components principle. Therefore, some of the sampling distributions addressed earlier (like the F distribution) probably do not accurately represent the distribution of the statistic they purport to be representing. In other words, they are only approximations to the true bivariate distribution and in violation of the condition $\alpha_1+\alpha_2=1$, required by the exact model (11).

Attempting to derive explicit expressions for the kth non-central moment of Z in approximation (23), as had been done for (1), failed to deliver an explicit expression.

### 3.2 A simplified approximate model – neglecting the normal component

If the error term in (8) is assumed small ($\varepsilon<<1$), most of the observed random variation can be assumed to originate in imperfect memory of the random phenomenon being modeled. We can then assume that (1) serves as good approximation to the true pdf of W, as the latter is given by (11) and (18). Identifying $\{\alpha,\lambda\}$ of (1) may be achieved by preserving the value of the pdf at the mode (eq. 3) and preserving the value of the mean (eq. 5). The normalizing parameter, $\kappa$, is determined by (6).

The resulting pdf is, similar to (1):

$$f(z;\kappa,\alpha,\lambda) = (\kappa)\exp\{-(\alpha)(\frac{1}{\lambda})[(z)^\lambda - 1]\},$$

$$\{-\frac{M}{\sigma} \leq z < \infty \; ; \; 0 \leq \lambda \leq UL\}$$

(30)

with $z=(w-M)/\sigma$, and $\{M,\sigma\}$ are the mode and the standard deviation of the modeled r.v., W. As shown earlier, pdfs of all distributions related to in this article, where the argument, w, appears only once, comply with this model. Examples are student's t and Cauchy (as shown earlier).

## 4. PREDICTION REGARDING A CERTIN COMPOUND DISTRIBUTION

A compound distribution is perhaps the best arena to examine the validity of the two-ignored-components paradigm. A general definition of a compound distribution is that it results from assuming that "a random variable is distributed according to some parametrized distribution, with the parameters of that distribution being assumed to be themselves random variables" (entry "Compound probability distribution" from Wikipedia, the free encyclopedia). The compound distribution is the result of marginalizing over the intermediate random variables that represent the parameters of the initial distribution (namely, integrating out the uncertain parameter(s)). Thus, the "identity" of the random variable being modeled is itself subject to random variation. This scenario roughly describes the reality described by the two-ignored-component principle, where it is assumed that lack of "inner" stability generates variation that is always part of the observed random variation (the term "roughly" is used since the theory underlying compound distributions in general does not relate to the roles of the exponential and the normal distributions, as elaborated on in this article; more on this further down).

An important type of compound distribution occurs when the parameter being marginalized over represents the number of random variables in a summation of random variables. When this number itself is a random variable we relate to a random sum, defined by:

$$S_N = X_1 + X_2 + .. + X_N \tag{31}$$

where the number of terms, N, is uncertain, $\{X_j\}$ are independent and identically distributed and each $X_j$ is independent of N. It is assumed that if N=0, we have $S_0=0$. By the law of total probability, the distribution function of S, the sum being independent of N, is a mixture distribution, with the mixture components being weighted by the probability function of N:

$$F_S(s) = \sum_{n=0}^{\infty} G_{Sn}(s) \Pr(N=n) = \sum_{n=0}^{\infty} G_{Sn}(s) h_N(n) \ , \tag{32}$$

where $F_S$ is the distribution function of S, $G_{S_n}$ is the (conditional) distribution function of $S_n$ (given that N=n) and $h_N$ is the probability function of N.

The mean and the variance of S can easily be shown to be:

$$E(S) = E(X)E(N)$$
$$Var(S) = E(N)Var(X) + Var(N)[E(X)]^2. \qquad (33)$$

Suppose now that $\{X_j\}$ are i.i.d exponentially distributed with parameter $\lambda$ and N is the discrete analog of the exponential distribution, namely, the geometric distribution with parameter p. According to the two-ignored-component principle, if p=1/2 then S is devoid of any memory (a scenario of lack of identity). Therefore we expect S to be exponential with mean equal to the standard deviation. From (33) we obtain:

$$E(S) = \frac{1-p}{\lambda p}$$
$$Var(S) = (\frac{1-p}{\lambda p})(\frac{1+p}{\lambda p}). \qquad (34)$$

As expected, for p=1/2 the mean is equal to the standard deviation and S is exponentially distributed. Note, that this result is expected under the new two-ignored-component paradigm. It could not have been predicted by any other existing general theory of random variation.

## 5. COMPOUND DISTRIBUTION AND THE TWO-COMPONENT PARADIGM

The concept of compound distribution seems to embody the two-component principle in the sense that the distribution, supposedly representing random variation associated with the modeled r.v., lacks "inner" stability due to its uncertain parameters. This may lead to the conclusion that the new paradigm is already addressed by current theory of compound distributions. It is therefore essential that we delineate the distinction between the two concepts (namely, "compound distribution" and "the two-component principle"):

**A.** The new paradigm claims that, apart from the exponential and the normal, *all* observed random variation, whether modeled by compound distribution or otherwise, contain two components indistinguishable from one another unless by statistical modeling. Furthermore, only the extreme scenarios of the exponential and the normal are true univariate distributions. All others are at least bivariate. The double appearance in the pdf of most univariate distributions of the distribution argument is remnant from the true bivariate distribution. This point has been elaborated on in subsection 3.1.

**B**. The new paradigm posts all current supposedly univariate distributions on a continuous spectrum that stretches from one "pure scenario" (exponential) to another (normal). The point on that spectrum occupied by a particular bivariate distribution is determined by how stable the underlying identity of the modeled random phenomenon is. For all cases, once the distribution ceases to be exponential, a normal component is generated, rendering the distribution truly bivariate with one r.v. representing the "identity" of the random phenomenon being modeled and the other representing the cumulative effect of external factors, exercised via a normal error (eq. 8). Current theory underlying compound distributions does not relate to the continuous spectrum from exponential to normal, implied by the two-ignored-components principle.

**C.** Some compound distributions (like those representing random sums) are in fact mixture distributions (refer to (32). The new paradigm delivers a bivariate distribution, where the two arguments of the joint pdf appear as components in a weighted average. However, the distribution itself is not a mixture.

# 6. FROM EXPONENTIAL TO NORMAL VIA GAMMA – AN ALTERNATIVE ROUTE TO THE NEW PARADIGM?

The gamma distribution is the distribution of the sum of N i.i.d exponential variables (N≥1), and it can be shown to asymptotically approximate the normal distribution. One may naturally ask whether the gamma distribution does not provide an alternative route to statistical modeling of the two-ignored-components principle? The answer is that being univariate, in violation of our basic supposition that **ALL** random variation (apart from that represented by the exponential and the normal) are bivariate, one cannot adopt the gamma, or better still the generalized gamma, but only as a univariate approximation to the true bivariate distribution.

## 7. CONCLUSION

The main claim of this article is that apart from random variation represented by the exponential and the normal – the only truly univariate distributions – all other seemingly univariate random variation is in fact bivariate. More specifically, the new paradigm perceives all random variation as emanating from two independent sources of random variation, internal and external, with two extreme scenarios having a single active source that deliver an exponential distribution (no externally generated random variation) or a normal distribution (no internally generated random variation).

A major concern that one may voice regarding the new paradigm is this:

Besides the logic underlying the new paradigm, which can be debated, what empirical evidence may one provide to corroborate the validity of the new paradigm?

As with all branches of science, a scientific theory cannot be proved but only be refuted. Since proving a theory is impossible, according to current philosophy of science, the only avenue to establish the validity of a new scientific theory is produce predictions and find out whether empirical observations are consistent with those predictions. In other words, compatibility between empirical observations and predictions, derived exclusively from a certain scientific theory, delivers the required evidence for the validity of the latter. We believe that ample empirical evidence has been provided in this article for the scientific validity of the new paradigm. It relies on five major empirical observations:

* The exponential and the normal have been shown to be special cases of a single RMM model; this affords considering all univariate distributions as truly bivariate distributions, positioned on a common continuous spectrum stretching from the exponential to the normal, with the bivariate distribution argument

appearing as a weighted average of two independent r.v.s (representing two independent sources of random variation) and with pdf given by (11);

* Numerous existent univariate distributions have been shown to be special cases of a univariate approximation to the bivariate true distribution (eq.23 and eq. 11, respectively);

* Appearance most often in current univariate distributions of the pdf argument twice indicates a bivariate source for the univariate distributions;

* For the generalized gamma and its four special cases (gamma, exponential, Weibull and chi-squared), as well as for Cauchy, the weights of approximation (23) sum up to one, namely: $\alpha_1+\alpha_2=1$, as predicted by (21); this result cannot be predicted unless by the two-component principle.

* Derivation of the exponential distribution as a special case of a random sum of N i.i.d exponential variables, with N having a geometric distribution with parameter p=1/2; this result is predicted from the two-component principle and could not have been predicted otherwise (unless via mathematical derivation, which then ceases to be a prediction from a general scientific theory).

Further empirical evidence may be assembled in the future.

Regarding the second observation above, perhaps it is useful to enumerate all distributions, derived directly or indirectly (via an approximation) from the two-component principle, using RMM models:

The normal, Pareto, Power function, uniform, generalized gamma (including special cases gamma, Weibull, exponential and chi-squared), F, log-normal, Student's t and Cauchy (thirteen distributions in all).

The main lesson from our ability to derive these distributions as approximate realizations of the new paradigm is that appearance of the r.v. argument in the expression for the pdf not once, as with the exponential and the normal, but twice (and never thrice) is residual reflection of the true nature of all observed random variation, namely, a two-component variation that results in a bivariate distribution. By this token, all current univariate distributions (besides the exponential and the normal) become simplifying distortion of the true nature of all observed random variation. Indeed, the sheer number of distributions having pdf that comply with (23) is the best corroborating evidence for the universal validity of the new paradigm for perceiving random variation.

The earliest we adopt this new insight, the more accurate and capable would our estimated models be in helping us predict random behavior and manipulate optimally scientific and managerial decisions in an environment of inherent uncertainty.